\documentclass[10pt]{amsart}
\usepackage{geometry}
\usepackage{graphicx}

\usepackage{tikz}
\usetikzlibrary{plotmarks}
\usepackage{pdfrender}
\usepackage{pifont}
\usepackage[hypcap]{caption}%ref of fiqure

\usepackage{float}
\newcommand*{\MarkSymbol}[3]{%
  \textpdfrender{%
    TextRenderingMode=FillStroke,
    LineWidth=.4pt,
    FillColor={#2},
    StrokeColor={#3},
  }{\scriptsize\ding{#1}}%
}
\newcommand*{\BlueRoute}{%
  \rotatebox{45}{%
    \MarkSymbol{110}{blue!60!white}{blue!80!black}%
  }%
}
\newcommand*{\WhiteSquare}{%
  \MarkSymbol{110}{white}{gray}%
}
\newcommand*{\RedStar}{%
  \MarkSymbol{72}{red!60!white}{red!80!black}%
}

\newcommand*{\BlackCirc}{%
  \MarkSymbol{108}{black}{black}%
}
\begin{document}

%%-----------------------------
%%      the top matter
%%-----------------------------http://link.springer.com/article/10.1007%2Fs00033-013-0380-7

\title[Numerical solutions for a Timoshenko-type system...]{Numerical solutions for a Timoshenko-type system with thermoelasticity with second sound}% At most 5 thanks
\author{Mohamed Ali Ayadi}\address{UR ANALYSE NON-LIN\'EAIRE ET G\'EOMETRIE, UR13ES32, Department of Mathematics, Faculty of Sciences of Tunis, University of Tunis El-Manar, 2092 El Manar II, Tunisia}
\email{ayadi.dali23@gmail.com}
\author{Ahmed Bchatnia}\address{UR ANALYSE NON-LIN\'EAIRE ET G\'EOMETRIE, UR13ES32, Department of Mathematics, Faculty of Sciences of Tunis, University of Tunis El-Manar, 2092 El Manar II, Tunisia}
\email{ahmed.bchatnia@fst.rnu.tn}
\author{Makram Hamouda}\address{Institute for Scientific Computing and Applied Mathematics, Indiana University, 831 E. 3rd St., Rawles Hall, Bloomington IN 47405, United States}
\email{mahamoud@indiana.edu}
%
%\dedicated{\it Dedicated to Maurice Dupont} %if necessary
%
\begin{abstract} In this work, we consider a nonlinear vibrating  Timoshenko system with thermoelasticity with second sound. We recall first the results of well-posdness and regularity  and the asymptotic behavior of the energy obtained in \cite{Ayadi}. Then, we use a fourth order finite difference scheme to compute the numerical solutions and thus we show the energy decay in several cases depending on the stability number.\\
{\bf{R\'esum\'e}} Dans ce travail, on consid\`ere le syst\`eme de Timoshenko non-lin\'eaire avec Thermo-\'elasticit\'e et deuxi\`eme son. On rappelle d'abord les r\'esultats d'existence, de r\'egularit\'e et du comportement asymptotique de l'\'energie obtenus dans \cite{Ayadi}. Ensuite, on valide num\'eriquement ces r\'esultats th\'eoriques. Pour cela, on utilise une m\'ethode de diff\'erences finies d'ordre $4$. Ainsi la solution num\'erique obtenue permet de valider la d\'ecroissance de l'\'energie dans plusieurs cas selon la valeur du param\`etre de stabilit\'e. 
\end{abstract}
\maketitle
%%-----------------------------
%%      your text
%%-----------------------------

\section{Introduction}
%\selectlanguage{english}

Historically, the first model of Timoshenko system was introduced in 1921 by Stephen Timoshenko in the absence of
dissipative term which describes the transverse vibration of the beam. Timoshenko considered thus the
following hyperbolic system:
\begin{equation}
\left\lbrace
\begin{array}{l}
\rho \varphi_{tt}=(k(\varphi_{x}+\psi ))_{x},\mbox{ }\,\hspace{3cm}\mbox{ in }(0,L)\times
\mathrm{I\hskip-2ptR}_{+}, \\
I_{\rho }\psi_{tt}=(EI\psi_{x})_{x}+k(\varphi_{t}+\psi ),\hspace{18mm}\mbox{ in }(0,L)\times \mathrm{I\hskip-2ptR}_{+},
\end{array}\right. \label{Timo1}
\end{equation}
where, $\rho $, $k$, $%
I_\rho $ and $EI$ are positive constants, $\varphi =\varphi (x,t)$ is the displacement vector and $\psi =\psi (x,t)$ is the rotation angle of the filament.

Among new works, many researchers used the classical model for the propagation of heat turns into the well-known equations for the temperature $\theta$ and the heat flux vector $q$
\begin{equation}
\theta_t +\beta div q=0, \label{chaleur}
\end{equation}
and
\begin{equation}
q+\kappa\nabla\theta =0.\label{fourier}
\end{equation}
with positive constants $\beta$ and $\kappa$. Substituting \eqref{fourier} (Fourier's low) into \eqref{chaleur}, yields the following parabolic heat equation
\begin{equation}
\theta_t -\beta \kappa ~\triangle q =0.
\end{equation}
Using the Fourier's low, Rivera and Racke \cite{MunozRiveraRackeMildy} investigated the following system
\begin{equation}
\left\lbrace
\begin{array}{l}
\rho _{1}\varphi _{tt}-k(\varphi _{x}+\psi )_{x}=0,\hspace{3cm}\mbox{ in }\,(0,L)\times \mathrm{I\hskip-2ptR}_{+}, \\
\rho _{2}\psi _{tt}-b\psi _{xx}+k(\varphi _{x}+\psi )+\delta \theta _{x}=0,\hspace{11mm}\mbox{ in }\,(0,L)\times \mathrm{I\hskip-2ptR}_{+}
,\\
\rho _{3}\theta _{t}-\kappa\beta\theta _{xx}+\delta \psi _{xt}=0,\hspace{2.88cm}\mbox{ in }(0,L)\times \mathrm{I\hskip-2ptR}_{+},
\end{array}\right.\label{four}
\end{equation}
where, $\rho_1 $, $\rho_2 $, $\rho_3 $, $k$, $
b $, $\kappa $, $\beta $, and $\delta $ are positive constants. They proved several exponential decay results for the linearized system and nonexponential stability results for the case of different wave speeds ($\frac{k}{\rho_1}\neq\frac{b}{\rho_2})$.\\
Later, Fern\'{a}ndez Sare and Racke considered in \cite{FernandezRacke} the following system:
\begin{equation}
\left\lbrace
\begin{array}{l}
\rho _{1}\varphi _{tt}-k(\varphi _{x}+\psi )_{x}=0,\hspace{5cm}\mbox{ in }\,
(0,L)\times  \mathrm{I\hskip-2ptR}_{+}, \\
\rho _{2}\psi _{tt}-b\psi _{xx}+k(\varphi _{x}+\psi )+\delta \theta
_{x}=0,\hspace{3,1cm}\mbox{ in }\, (0,L)\times  \mathrm{I\hskip-2ptR}_{+}, \\
\rho _{3}\theta _{t}+ q_{x}+\delta \psi _{tx}=0,\hspace{5,45cm}\mbox{ in }\,
(0,L)\times  \mathrm{I\hskip-2ptR}_{+}, \\
\tau q_{t}+\beta q+ \theta _{x}=0,\hspace{5,9cm}\mbox{ in }\,  (0,L)\times  \mathrm{I\hskip-2ptR}_{+},\\
 \varphi(0,t)\!\!=\!\!\varphi(L,t)\!\!=\!\!\psi_x(0,t)\!\!=\!\!\psi_x(L,t)
\!\!=\!\!\theta_x(0,t)=\!\! \theta_x(L,t)=0,\hspace{3mm} \forall\ t\in  \mathrm{I\hskip-2ptR}_{+},
\end{array}\right. \label{Timo6}
\end{equation}
and proved that the coupling via Cattaneo's law $(\ref{Timo6})_4$ does not make the energy decays exponentially which is usually obtained for the coupling via Fourier's law (system (\ref{four})).

Numerically, Raposo et al \cite{Raposo} considered the following  Timoshenko system with a delay term in the feedback:
\begin{equation}
\left\lbrace
\begin{array}{l}
\rho_1 \varphi_{tt}(x,t)-k(\varphi_{x}+\psi ))_{x}(x,t)+\mu_1\varphi_t(x,t)+\mu_2\varphi_t(x,t-\tau)=0,\mbox{ }\,\hspace{1,35cm}\mbox{ in }(0,L)\times
\mathrm{I\hskip-2ptR}_{+}, \\
\rho_{2}\psi _{tt}(x,t)-b\psi _{xx}(x,t)-k(\varphi_{t}+\psi )(x,t)+\mu_3\psi_t(x,t)+\mu_4\psi_t(x,t-\tau)=0,\mbox{ in
}(0,L)\times \mathrm{I\hskip-2ptR}_{+},\\
 \varphi(0,t)=\varphi(L,t)=\psi(0,t)=\psi(L,t)=0,\hspace{5,35cm}\forall \ t>0,
\end{array}\right.
\end{equation}
and they gave different tests of decay results for the solutions of the previous system.

Recently, Ayadi et al. \cite{Ayadi} considered the following coupling of two wave equations of Timoshenko
type system:
\begin{equation}
\left\{
\begin{array}{l}
\rho _{1}\varphi _{tt}-k(\varphi _{x}+\psi )_{x}=0,\hspace{4.3cm}\textnormal{in }%
(0,1)\times \mathrm{I\hskip-2ptR}_{+}, \\
\rho _{2}\psi _{tt}-b \psi _{xx}+k(\varphi _{x}+\psi )
+\delta \theta _{x}+\alpha (t)h(\psi _{t})=0,\hspace{0.5cm} \textnormal{in }(0,1)\times
\mathrm{I\hskip-2ptR}_{+}, \\
\rho _{3}\theta _{t}+q_{x}+\delta \psi _{xt}=0,\hspace{4.72cm} \textnormal{in }(0,1)\times
\mathrm{I\hskip-2ptR}_{+}, \\
\tau q_{t}+\beta q+\theta _{x}=0,\hspace{5.18cm} \textnormal{in }(0,1)\times \mathrm{I\hskip%
-2ptR}_{+}.
\end{array}
\right.
\label{2}
\end{equation}
In order to study the stability properties of the solution of the system (\ref{2}), the authors introduced in \cite{Ayadi} a
stability number
\[
\mu =\left( \tau -\frac{\rho _{1}}{k\rho _{3}}\right) \left( \frac{\rho _{2}}{b}-
\frac{\rho _{1}}{k}\right) -\frac{\tau\delta ^{2}\rho _{1}}{bk\rho _{3}}.
\]
This number $\mu$ is crucial in determining the asymptotic behavior of the energy associated with system (\ref{2}).

This paper is organized as follows: in Section 2 we recall the results of the existence and asymptotic behavior of the solutions of the system (\ref{2}). In Section 3 we present the numerical solutions in some particular cases.
\section{Results of existence and asymptotic behavior}
In this section, we recall the results obtained in \cite{Ayadi}. Precisely, the authors studied the system:
\begin{equation}
\left\{
\begin{array}{l}
\rho _{1}\varphi _{tt}-k(\varphi _{x}+\psi )_{x}=0,\hspace{4.3cm}\textnormal{in }%
(0,1)\times \mathrm{I\hskip-2ptR}_{+}, \\
\rho _{2}\psi _{tt}-b \psi _{xx}+k(\varphi _{x}+\psi )
+\delta \theta _{x}+\alpha (t)h(\psi _{t})=0,\hspace{0.5cm} \textnormal{in }(0,1)\times
\mathrm{I\hskip-2ptR}_{+}, \\
\rho _{3}\theta _{t}+q_{x}+\delta \psi _{xt}=0,\hspace{4.718cm} \textnormal{in }(0,1)\times
\mathrm{I\hskip-2ptR}_{+},\hspace{0.3cm} \\
\tau q_{t}+\beta q+\theta _{x}=0,\hspace{5.18cm} \textnormal{in }(0,1)\times \mathrm{I\hskip%
-2ptR}_{+},
\end{array}
\right.
\label{1}
\end{equation}
where, $\rho _{1}$, $%
\rho _{2}$, $\rho _{3}$, $b$, $k$, $\delta $, $\beta $ are positive constants, $\varphi =\varphi (x,t)$ is the displacement vector, $\psi =\psi (x,t)$
is the rotation angle of the filament, $\theta =\theta (x,t)$ is the
temperature difference and $q=q(x,t)$ is the heat flux vector.
Also, $ \alpha $ and $h$ verify the assumptions:

$(A_{1})$  $\alpha $: $\mathbb{R}_{+}\rightarrow \mathbb{R}_{+}$ is a differentiable and decreasing function.

$(A_{2})$  $h$: $\mathbb{R}\rightarrow \mathbb{R}$ is a continuous non-decreasing function with $h(0)=0$ and there exists a continuous strictly increasing odd function $h_{0}\in C([0,+\infty))$, continuously differentiable in a neighborhood of $0$, satisfying $h_{0}(0)=0$ and such that
\newline
$$ \left\{
\begin{array}{l}
h_{0}( s )\leq \vert h(s) \vert \leq h_{0}^{-1}(s ),%
\hspace{0.65cm} \text{for all}\hspace{0.2em} \vert s\vert\leq\varepsilon , \\
c_{1} \vert s \vert\leq \vert h(s)\vert \leq c_{2} \vert s\vert ,\hspace{%
2.8em}\text{ for all} \hspace{0.2em}\vert s\vert\geq\varepsilon,%
\end{array}
\right. $$ \\
where $c_{i} > 0$ for i = 1, 2. \\
With \eqref{1}, we associate the boundary conditions given by
\begin{equation}
\varphi _{x}(0,t)=\varphi _{x}(1,t)=\psi (0,t)=\psi (1,t)=q(0,t)=q(1,t)=0,%
\hspace{0.5cm}\forall\mbox{ }t\geq 0 ,
\label{1b}
\end{equation}
and the following initial conditions
\begin{equation}
\left\{
\begin{array}{l}
\varphi (x,0)=\varphi _{0}(x),\mbox{ }\varphi _{t}(x,0)=\varphi _{1}(x),\hspace{0.8cm}
\forall\mbox{ }x\in (0,1) ,\\
\psi (x,0)=\psi _{0}(x),\mbox{ }\psi _{t}(x,0)=\psi _{1}(x),\hspace{0.8cm}\forall\mbox{ }x\in
(0,1), \\
\theta (x,0)=\theta _{0}(x),\mbox{ }q(x,0)=q_{0}(x),\hspace{1.2cm}\forall\mbox{ }x\in (0,1).
\end{array}
\right.
\label{1i}
\end{equation}
\subsection{Well-posedness and regularity}
Here, we state the existence and uniqueness results of solutions of the Timoshenko system composed of (\ref{1}), (\ref{1b}) and (\ref{1i}) (see \cite{Ayadi}).

\textbf{Theorem 1}.
\textit{Assume that $(A_{1})$ and $(A_{2})$ are satisfied, then for all initial data \\ \mbox{} \hspace{2cm} $(\varphi_{0},\varphi_{1},\psi_{0},\psi_{1},\theta_{0},q_{0})\in (H^{2}_{\star}(0,1)\cap H_{\star}^{1}(0,1))\times H_{\star}^{1}(0,1)\times (H^{2}(0,1)\cap H_{0}^{1}(0,1))$ \\
\mbox{} \hspace{5,6cm} $\times H_{0}^{1}(0,1)\times H_{\star}^{1}(0,1)\times H_{0}^{1}(0,1),$\\
the system (\ref {1})--(\ref{1i}) has a unique solution $(\varphi,\psi,\theta,q)$ that verifies  \\
$\left. \right. \hspace{2cm}
(\varphi,\psi)\in C^{0}(\mathbb{R}_{+},(H^{2}_{\star}(0,1)\cap H_{\star}^{1}(0,1))\times(H^{2}(0,1)\cap H_{0}^{1}(0,1)))$ \\
$ \left. \right. \hspace{3cm}\cap \ C^{1}(\mathbb{R}_{+},H_{\star}^{1}(0,1)\times H_{0}^{1}(0,1))\cap \ C^{2}(\mathbb{R}_{+},L^{2}_{\star }(0,1)\times L^{2}(0,1)),$\\
and $$\hspace{1.8cm}(\theta,q)\in  \ C^{0}(\mathbb{R}_{+},H_{\star}^{1}(0,1)\times H_{0}^{1}(0,1))\cap \ C^{1}(\mathbb{R}_{+},L^{2}_{\star }(0,1)\times L^{2}(0,1)),$$
where
$$L_\star^2 (0,1)=\{v\in L^2 (0,1)  \mbox{ s.t. }  \int_{0}^{1}v(s)ds=0 \},$$
$$H_\star^1 (0,1)= H^1 (0,1)\cap L_\star ^2(0,1),$$ and
$$H_\star^2 (0,1)=\{v\in H^2 (0,1) \mbox{ s.t. } v_x(0)=v_x(1)=0 \}.$$}
\subsection{Asymptotic behavior }
In this subsection, we give the general decay results for a wide class of relaxation functions (denoted here by $h$).

\textbf{Theorem 2}. \textit{Let us suppose that $(A_{1})$ and $(A_{2})$ are satisfied, then for
$\mu =0$ there exist positive constants $k_{1}$, $k_{2}$, $k_{3}$ and $\varepsilon_{0}$ such that the energy $E(t),$ associated with \eqref{1}--\eqref{1i}, satisfies
\begin{equation}
\left. E(t)\leq k_{3} H^{-1}_{1} \left( k_{1} \int_{0}^{t} \alpha(s)\ ds+
k_{2} \right), \hspace{1.cm} \text{for all} \ \ t\geq 0, \right.  \label{theoremenergy1}
\end{equation}
where
\[
H_{1}(t) = \int_{t}^{1} \frac{1}{H_{2}(s)}ds,\hspace{1cm} H_{2}(t)=tH^{^{%
\prime }}(\varepsilon_{0}t).
\]
Here $H_{1}$ is a strictly decreasing and convex function on $(0,1]$, with $\displaystyle\lim_{t\rightarrow0}H_{1}(t)= +\infty$.}

In the following, in order to show explicit stability results in term of asymptotic profiles in time, we consider some special values for the function $h_0$.

\begin{itemize}

\item[\underline{Example 1.}]

For $h_{0}(s)=cs^{p} $, we have\\
$\bullet$ If $p=1$, then $E(t)\leq k_{3}\exp(-c(k_{1}\int^{t}_{0}\alpha(s)\
ds+k_{2})).$\\
$\bullet$ If $p>1$, then $E(t)\leq c( k_{1}\int_{0}^{t}\alpha(s)\ ds +k_{2})^{-\frac{2}{p-1}}.$

\item[\underline{Example 2.}] For $h_{0}(s)=\exp(-\frac{1}{s})$, we have\\
 $E(t)\leq k_{3}\varepsilon_{0}^{-1}\left(  \ln\left( \frac{k_{1} \int_{0}^{t}\alpha(s)ds +k_{2}+c\exp(\frac{1}{\sqrt{\varepsilon_{0}}})}{c}\right) \right)^{-2}.$%
\item[\underline{Example 3.}]
For $h_{0}(s)=\frac{1}{s}\exp(-\frac{1}{s^{2}})$, we have \\
$E(t)\leq \varepsilon \left( \ln(\frac{k_{1}\int_{0}^{t}\alpha(s)ds +k_{2}+c\exp(\frac{1}{\varepsilon_{0}})}{c})\right)^{-1}.$

\end{itemize}

Next, we will consider the case where the stability number $\mu \neq 0.$

\textbf{Theorem 3}.\label{th3}
\textit{Let us suppose that the derivative of the function $h$ is bounded and the assumptions $(A_{1})$ and $(A_{2})$ hold, then for $\mu \neq 0,$ the energy solution of \eqref{1}--\eqref{1i} satisfies}
\begin{equation}
 E(t)\leq E(0)H_{2}^{-1}(\frac{c}{t}),\label{theoremenergy2}
\end{equation}
\textit{ where
$$
H_{2}(t)=tH^{^{\prime }}(\varepsilon_{0}t) \mbox{ with } \displaystyle\lim_{t\rightarrow0}H_{2}(t)=0.
$$ }
In the following, we give some examples to illustrate the energy decay rates
given by Theorem 3.

\underline{Example 1.}
For $h_{0}(s)=cs^{p} $, then \\
$\bullet$ If p=1, we have $E(t)\leq \frac{c}{t}$.\\
$\bullet$ If $p>1$ we have $
E(t)\leq ct^{-\frac{2}{p+1}}.
$

\underline{Example 2.}
 Let $h$ be given by $h(x)=\frac{1}{x^{3}}\exp(-\frac{1}{x^{2}})$ and we choose $h_{0}(x)=\frac{1+x^2}{x^3}\exp(-\frac{1}{x^{2}})$, we obtain
 \begin{center}
 $E(t)\leq c(\ln(t))^{-1}.$
 \end{center}

\section{Numerical solution}
We will start making use of Finite Difference Method to derive a discrete representation of the solution of the Timoshenko system \eqref{1}--\eqref{1i} in the particular case $\alpha(t)=1$ and $h(s)=s$. \\ More precisely, we use the classical finite difference discretization for the  temporal variable and the Implicit Compact Finite Difference Method of fourth-order for discretization of the space variable. The full nonlinear case and the comparison between the different Difference Finite Methods (implicit, explicit and semi-implicit) will be considered in a subsequent work.
\subsection{Discrete formulation}
Consider the discrete domain of $\Omega_{h}=(0,1)$ with uniform grid $ x_{i}=ih,\  i=0,1,...,I; h=\frac{1}{I} $. The temporal discretization of the interval $T_n =(0,T)$ is given by $ t_{n}=n\kappa,\ n=0,1,...,N;\ \kappa=ch$, where $c$ is a positive constant and $I$ and $N$ are two positive integers. Denote by $\omega (x_i,t_n)=\omega_{i}^{n}$ the value of the function $\omega$ evaluated at the point $x_i$ and the instant $t_n$.\\
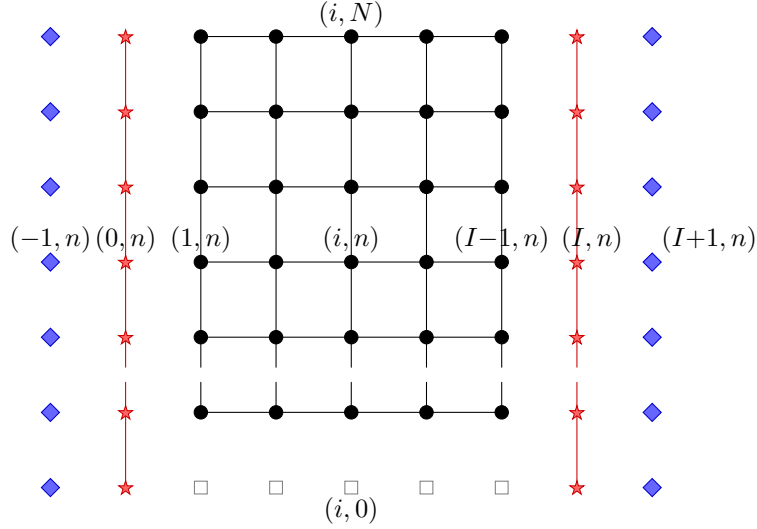
\begin{figure}[H]
\begin{tikzpicture}[line cap=round]
  % Lines
  \draw[red!80!black]
    (-3, 0) -- (-3, 1.4) (-3, 1.6) -- (-3, 6)
    (3, 0) -- (3, 1.4) (3, 1.6) -- (3, 6)
  ;
  \draw[black]
    \foreach \x in {-2, ..., 2} {
      (\x, 1) -- (\x, 1.4) (\x, 1.6) -- (\x, 6)
    }
    \foreach \y in {1, ..., 6} {
      (-2, \y) -- (2, \y)
    }
  ;
  \path
    \foreach \x in {-4, 4} {
      \foreach \y in {0, ..., 6} {
        (\x, \y) node{\BlueRoute}
      }
    }
    \foreach \x in {-3, 3} {
      \foreach \y in {0, ..., 6} {
       (\x, \y) node{\RedStar}
      }
    }
    \foreach \x in {-2, ..., 2} {
       (\x, 0) node {\WhiteSquare}
       \foreach \y in {1, ..., 6} {
         (\x, \y) node {\BlackCirc}
       }
    }
  ;
  \node[below] at (0, 0) {$(i, 0)$};
  \node[above] at (0, 6) {$(i, N)$};
  \node[above] at (-4, 3) {$(-1, n)$};
  \node[above] at (-3, 3) {$(0, n)$};
  \node[above] at (-2, 3) {$(1, n)$};
  \node[above] at (0, 3) {$(i, n)$};
  \node[above] at (2, 3) {$(I{-}1, n)$};
  \node[above, xshift=2mm] at (3, 3) {$(I, n)$};
  \node[above right] at (4, 3) {$(I{+}1, n)$};
\end{tikzpicture}
\caption{\ Mesh of the domain $\Omega_{h} \times T_{n} $}
\end{figure}
In Figure 1, we show the pattern mesh of $\omega$ using the discretization of the intervals $(0,1)$ and $(0,T),$ with the classification of nodes is as follows: internal (circles), boundaries (stars), initials (squares) and ghost (diamonds).\\
Now, We define the following approximation of the derivatives of $\omega$
\begin{eqnarray}
(\omega_{t})_{i}^{n}\simeq \frac{\omega_{i}^{n+1}-\omega_{i}^{n-1}}{2\kappa},\hspace{0.5cm}(\omega_{tt})_{i}^{n}\simeq \frac{\omega_{i}^{n+1}-2\omega_{i}^{n}+\omega_{i}^{n-1}}{\kappa^2},\label{discri1}
\end{eqnarray}
\begin{eqnarray}
(\omega_{x})_{i}^{n}\simeq\frac{\omega_{i+1}^{n}-\omega_{i-1}^{n}}{2h},\hspace{0.7cm}(\omega_{xx})_{i}^{n}\simeq\frac{1}{h^2}\left[ \frac{\delta_{x}^{2}}{1+\dfrac{1}{12}\delta_{x}^{2}}\right]\omega_{i}^{n}, \hspace{1.8cm}\label{discri2}
\end{eqnarray}
with
\begin{eqnarray}
 \hspace{0.3cm} \delta_{x}^{2}\omega_{i}^{n}=\omega_{i+1}^{n}-2\omega_{i}^{n}+\omega_{i-1}^{n} \text{ and}\left[ 1+\dfrac{1}{12}\delta_{x}^{2}\right]\omega_{i}^{n}=\frac{1}{12}\omega_{i+1}^{n}+\frac{5}{6}\omega_{i}^{n}+\frac{1}{12}\omega_{i-1}^{n}. \hspace{0.5cm}
\label{discri3}
\end{eqnarray}
Using \eqref{discri1} and \eqref{discri2}, we obtain the discrete formulation of the system \eqref{1} as follows:
\begin{equation}
\left\lbrace
\begin{array}{l}
\rho_1\frac{\varphi_i^{n+1}-2\varphi_i^{n}+\varphi_i^{n-1}}{\kappa^2}-\frac{k}{h^2}\left[ \frac{\delta_x^2}{1+\frac{1}{12}\delta_x^2}\right]\varphi_{i}^{n}-\frac{k}{2h}(\psi_{i+1}^{n}-\psi_{i-1}^{n}) =0, \vspace{0.2cm}\\

\rho_2\frac{\psi_i^{n+1}-2\psi_i^{n}+\psi_i^{n-1}}{\kappa^2}-\frac{b}{h^2}\left[ \frac{\delta_x^2}{1+\frac{1}{12}\delta_x^2}\right]\psi_{i}^{n}+\frac{k}{2h}(\varphi_{i+1}^{n}-\varphi_{i-1}^{n})+k\psi_i^n \vspace{0.2cm}\\
\hspace{2cm}+\frac{\delta}{2h}(\theta_{i+1}^{n}-\theta_{i-1}^{n})+\frac{1}{2\kappa}(\psi_{i}^{n+1}-\psi_{i}^{n-1}) =0,\vspace{0.2cm}\\

\rho_3\frac{(\theta_{i}^{n+1}-\theta_{i}^{n-1})}{2\kappa}+\frac{1}{2h}(q_{i+1}^{n}-q_{i-1}^{n})+\frac{\delta}{4\kappa h}(\psi_{i+1}^{n+1}-\psi_{i-1}^{n+1}) -\frac{\delta}{4\kappa h}(\psi_{i+1}^{n-1}-\psi_{i-1}^{n-1})=0,\vspace{0.2cm}\\

\tau\frac{(q_{i}^{n+1}-q_{i}^{n-1})}{2\kappa}+\beta q_i^n+\frac{1}{2h}(\theta_{i+1}^{n}-\theta_{i-1}^{n})=0.
\end{array}\right. \label{syseq}
\end{equation}
Multiplying $\eqref{syseq}_1$ and $\eqref{syseq}_2$ by $\kappa^2[1+\frac{1}{12}\delta_x^2]$ and using \eqref{discri3}, we obtain
\begin{equation}
\left\lbrace
\begin{array}{l}
\rho_1(\frac{1}{12}\varphi_{i+1}^{n+1}+\frac{5}{6}\varphi_i^{n+1}+\frac{1}{12}\varphi_{i-1}^{n+1})-\frac{1}{6}(\rho_1+6a_1)\varphi_{i+1}^{n}-\frac{1}{3}(5\rho_1-6a_1)\varphi_{i}^{n}-\frac{1}{6}(\rho_1+6a_1)\varphi_{i-1}^{n}\vspace{0.2cm}\\ 
+\rho_1(\frac{1}{12}\varphi_{i+1}^{n-1}+\frac{5}{6}\varphi_i^{n-1}+\frac{1}{12}\varphi_{i-1}^{n-1})-a_2(\frac{1}{12}\psi_{i+2}^{n}+\frac{5}{6}\psi_{i+1}^{n}-\frac{5}{6}\psi_{i-1}^{n}-\frac{1}{12}\psi_{i-2}^{n})=0,\vspace{0.6cm}\\ 
(\rho_2+b_3)(\frac{1}{12}\psi_{i+1}^{n+1}+\frac{5}{6}\psi_i^{n+1}+\frac{1}{12}\psi_{i-1}^{n+1})-\frac{1}{12}(2\rho_2+12b_1-b_0)\psi_{i+1}^{n}-\frac{1}{6}(10\rho_2-12b_1-5b_0)\psi_{i}^{n}\vspace{0.2cm}\\
-\frac{1}{12}(2\rho_2\!+\!12b_1\!-\!b_0)\psi_{i-1}^{n}\!+\!b_2(\frac{1}{12}\varphi_{i+2}^{n}+\frac{5}{6}\varphi_{i+1}^{n}-\frac{5}{6}\varphi_{i-1}^{n}\!-\!\frac{1}{12}\varphi_{i-2}^{n})\vspace{0.2cm}\\ +b_4(\frac{1}{12}\theta_{i+2}^{n}+\frac{5}{6}\theta_{i+1}^{n}-\frac{5}{6}\theta_{i-1}^{n}-\frac{1}{12}\theta_{i-2}^{n})-(-\rho_2+b_3)(\frac{1}{12}\psi_{i+1}^{n-1}+\frac{5}{6}\psi_i^{n-1}+\frac{1}{12}\psi_{i-1}^{n-1})=0,
\label{discri4}
\end{array}\right. 
\end{equation}
where the parameters are defined by
$$a_1=\frac{k\kappa^2}{h^2},\hspace{0.5cm}a_2=b_2=\frac{k\kappa^2}{2h},\hspace{0.5cm}b_0=k\kappa^2,\hspace{0.5cm}b_1=\frac{b\kappa^2}{h^2},\hspace{0.5cm}b_3=\frac{\kappa}{2}\hspace{0.4cm}\text{and }\hspace{0.4cm}b_4=\frac{\delta\kappa^2}{2h}.$$
The discrete formulation of the initial conditions \eqref{1i} is defined by
\begin{equation*}
\left\lbrace \begin{array}{c}
\varphi(x_{i},0)=\varphi_{i}^{0}=(\varphi_{0})_{i},\hspace{0.3cm}\psi(x_{i},0)=\psi_{i}^{0}
=(\psi_{0})_{i}, \hspace{1cm}\text{for all } x_{i}\in  \mathring{\Omega_h},\\
\varphi_{t}(x_{i},0)=(\varphi_{t})_{i}^{0}
 =(\varphi_{1})_{i},
 \hspace{0.2cm}\psi_{t}(x_{i},0)=(\psi_{t})_{i}^{0}=(\psi_{1})_{i} ,\hspace{0.2cm}\text{for all }  x_{i}\in  \mathring{\Omega_h}.\end{array}\right.
 \end{equation*}
 The discrete formulation of the boundary conditions \eqref{1b} defined by
\begin{equation*} (\varphi_{x})_{0}^{n}=(\varphi_{x})_{I}^{n}=\psi_{0}^{n}=\psi_{I}^{n}=q_{0}^{n}=q_{I}^{n}=0,\hspace{0.3cm}\text{for all } t_{n}\in T_n . \label{discri6}
 \end{equation*}
In addition, it is natural to assume that $\varphi_{-1}^n=\varphi_0^n=\varphi_1^n$ and $\theta_{-1}^n=\theta_0^n=\theta_1^n$, since we have  $(\theta_{x})_{0}^{n}=(\theta_{x})_{I}^{n}=0$, thanks to $\eqref{1}_4$. Hence, we obtain the following linear algebraic system:
\begin{eqnarray}
\left\{
\begin{array}{c}
A_{1} \Phi^{n+1}=B_1 \Phi^{n}+C_1 \Psi^{n}+D_1 \Phi^{n-1},\\
\hspace{1.3cm} A_{2} \Psi^{n+1}=B_2 \Psi^{n}+C_2 \Phi^{n}+D_2\Psi^{n-1}+F_{2}\Theta^{n},\\
\hspace{0cm}A_{3} \Theta^{n+1}+L_{3}\Psi^{n+1}=B_3 \Theta^{n-1}-C_3 Q^{n}+D_3 \Psi^{n-1},\\
\hspace{0.2cm}A_{4} Q^{n+1}=B_4 Q^{n-1}-C_4 Q^{n}-D_4 \Theta^{n},\\

 \end{array}%
\right.\label{systemediscri}
\end{eqnarray}
with $ \Phi^{n}=(\varphi_{1}^{n},\varphi_{2}^{n},...,\varphi_{I-1}^{n})^{t}$, $\Psi^{n}=(\psi_{1}^{n},\psi_{2}^{n},...,\psi_{I-1}^{n})^{t}$, $ \Theta^{n}=(\theta_{1}^{n},\theta_{2}^{n},...,\theta_{I-1}^{n})^{t}$, $Q^{n}=(q_{1}^{n},q_{2}^{n},...,q_{I-1}^{n})^{t}$,  $\text{for all } n \in \{0,1...,N-1\}$ and $A_p$, $B_p$, $C_p$, $D_p$, $F_2$ and $L_3$ are $(I-1)$ square matrices for $p=1,...,4$ which will be defined below. First, we have
$$A_1=\begin{pmatrix}
 \frac{11}{12}\rho_1 & \frac{\rho_1}{12} & 0 & \cdots & \cdots & \cdots & 0 \\
                   \frac{\rho_1}{12} & \frac{5}{6}\rho_1 & \ddots & \ddots &  &  & \vdots \\
                   0 & \ddots & \ddots & \ddots & \ddots  &  & \vdots \\
                   \vdots & \ddots & \ddots & \ddots & \ddots & \ddots  & \vdots \\
                   \vdots & & \ddots & \ddots & \ddots & \ddots  & 0 \\
                   \vdots & & & \ddots & \ddots & \frac{5}{6}\rho_1 & \frac{\rho_1}{12} \\
                  0 & \cdots & \cdots & \cdots & 0 & \frac{\rho_1}{12} & \frac{11}{12}\rho_1
\end{pmatrix}=-D_1.$$
It is clear that $A_1$ is almost tridiagonal matrix except in the first and last diagonal terms where the coefficient is  $\frac{11}{12}\rho_1.$ Similarly, the matrix $B_1$ is an almost tridiagonal matrix given as follows:
$$B_1=\begin{pmatrix}
 \beta_3 & \beta_1 & 0 & \cdots & \cdots & \cdots & 0 \\
                   \beta_1 & \beta_2 & \beta_1 & \ddots &  &  & \vdots \\
                   0 & \ddots & \ddots & \ddots & \ddots  &  & \vdots \\
                   \vdots & \ddots & \ddots & \ddots & \ddots & \ddots  & \vdots \\
                   \vdots & & \ddots & \ddots & \ddots & \ddots  & 0 \\
                   \vdots & & & \ddots & \ddots & \beta_2 & \beta_1 \\
                  0 & \cdots & \cdots & \cdots & 0 & \beta_1 & \beta_3
\end{pmatrix},$$
where  $\beta_1=\frac{1}{6}(\rho_1+6a_1)$, $\beta_2=\frac{1}{3}(5\rho_1-6a_1)$ and $\beta_3=\frac{1}{6}(11\rho_1-6a_1)$.\\
Now, the matrices $C_1$, $A_2$, $B_2$ and $D_2$ are given by: \vspace{-1mm}
$$C_{1}=pentadiag(-\frac{1}{12}a_2,-\frac{5}{6}a_2,0,\frac{5}{6}a_2,\frac{1}{12}a_2),\vspace{-1mm}$$
$$A_{2}=tridiag(\frac{1}{12}(\rho_2+b_3),\frac{5}{6}(\rho_2+b_3),\frac{1}{12}(\rho_2+b_3)),\vspace{-1mm}$$
$$B_{2}=tridiag(\frac{1}{12}(2\rho_2+12b_1-b_0),\frac{1}{6}(10\rho_2-12b_1-5b_0),\frac{1}{12}(2\rho_2+12b_1-b_0)),\vspace{-1mm}$$
$$D_{2}=tridiag(\frac{1}{12}(-\rho_2+b_3),\frac{5}{6}(-\rho_2+b_3),\frac{1}{12}(-\rho_2+b_3)).$$
However, the matrices $C_2$ and $F_2$ do not have any particular form and they are given as follows:
$$C_2=\begin{pmatrix}
  \frac{11b_2}{12} & -\frac{5b_2}{6} & -\frac{b_2}{12} & 0 & \cdots & \cdots & 0  \\
                 \frac{11b_2}{12} & 0 & -\frac{5b_2}{6} & -\frac{b_2}{12} & \ddots &  & \vdots \\
                   \frac{b_2}{12}  & \frac{5b_2}{6} & \ddots & \ddots & \ddots  & \ddots & \vdots \\
                   0 & \ddots & \ddots & \ddots & \ddots & \ddots  & 0 \\
                   \vdots & \ddots & \ddots & \ddots & \ddots & -\frac{5b_2}{6}  & -\frac{b_2}{12}  \\
                   \vdots & & \ddots & \ddots & \ddots & 0 & -\frac{11b_2}{12}  \\
                  0 & \cdots & \cdots & 0 & \frac{b_2}{12}  &  \frac{5b_2}{6} & -\frac{11b_2}{12}
\end{pmatrix},$$
$$F_2=\begin{pmatrix}
  \frac{11b_4}{12} & -\frac{5b_4}{6} & -\frac{b_4}{12} & 0 & \cdots & \cdots & 0  \\
                 \frac{11b_4}{12} & 0 & -\frac{5b_4}{6} & -\frac{b_4}{12} & \ddots &  & \vdots \\
                   \frac{b_4}{12}  & \frac{5b_4}{6} & \ddots & \ddots & \ddots  & \ddots & \vdots \\
                   0 & \ddots & \ddots & \ddots & \ddots & \ddots  & 0 \\
                   \vdots & \ddots & \ddots & \ddots & \ddots & -\frac{5b_4}{6}  & -\frac{b_4}{12}  \\
                   \vdots & & \ddots & \ddots & \ddots & 0 & -\frac{11b_4}{12}  \\
                  0 & \cdots & \cdots & 0 & \frac{b_4}{12}  &  \frac{5b_4}{6} & -\frac{11b_4}{12}
\end{pmatrix}.$$
Finally, we denote by $I_d$ the identity matrix of size $(I-1)$ and we introduce the following parameters
$$\tau_1=\frac{\rho_3}{2\kappa},\hspace{1cm}\tau_2=\frac{1}{2h},\hspace{1cm}\tau_3=\frac{\delta}{4\kappa h},\hspace{1cm}r_1=\frac{\tau}{2\kappa},\hspace{1cm}r_2=\beta \mbox{  \hspace{0.2cm}   and } \hspace{0.7cm}r_3=\frac{1}{2h}.$$
Finally, we define the remaining matrices as below:
$$A_3=B_3=\tau_1I_d,$$
$$L_3=D_3=tridiag(-\tau_3,0,\tau_3),$$
$$C_3=tridiag(-\tau_2,0,\tau_2),$$
$$A_4=B_4=r_1I_d,$$
$$C_4=r_2I_d,$$
$$D_4=\begin{pmatrix}
 -r_3 & r_3 & 0 & \cdots & \cdots & \cdots & 0 \\
                  - r_3 & 0 & r_3 & \ddots &  &  & \vdots \\
                   0 & \ddots & \ddots & \ddots & \ddots  &  & \vdots \\
                   \vdots & \ddots & \ddots & \ddots & \ddots & \ddots  & \vdots \\
                   \vdots & & \ddots & \ddots & \ddots & \ddots  & 0 \\
                   \vdots & & & \ddots & \ddots & 0 & r_3 \\
                  0 & \cdots & \cdots & \cdots & 0 & -r_3 & r_3
\end{pmatrix}.$$
\subsection{Numerical tests}
To verify the asymptotic behavior of the solutions of the Timoshenko system \eqref{1}, we consider  the following data $I =26$, $T=35$,  $c=0,05$ and the initial conditions:\\
\begin{eqnarray}
\left\lbrace \begin{array}{l}
 \varphi_{0}(x)=\psi_{0}(x)= \theta_{0}(x)=q_{0}(x)=0,\\
 \varphi_{1}(x)=cos(\pi x),\hspace{0.3cm} \psi_{1}(x)=sin(2\pi x),\\
 \theta_t(x,0)=\theta_{1}(x)=(\frac{-2\pi\delta }{\rho_3})cos(2\pi x),\\
 q_t(x,0)=q_{1}(x)=0.
  \end{array}\right.
\end{eqnarray}
Note that in what follows the energy decay of the solution is proven by taking the maximum value of the function of the displacement $\varphi$.
\subsubsection{The case $\mu = 0$.}
For the following numerical computation we will consider different values for the parameters $k,\rho_{1},\rho_{2},\rho_{3},b,\beta $ and $\tau $. For example, in Figure 2 below, we take $k=\rho_{1}=\rho_{2}=2$,\ $b=\rho_{3}=\beta =1$,\ $\delta =\sqrt{\frac{2}{3}}$ and $\tau  =3$.
 \begin{figure}[H]
        \centering
        \includegraphics[width=4.2cm]{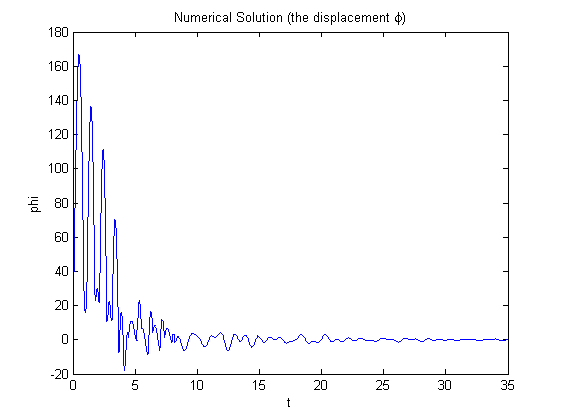}
        \caption{\label{fig1}}
\end{figure}
We recall here that we have theoretically obtained an exponential decay of the energy of the Timoshenko solution (see Theorem 2).
\subsubsection{The case $\mu \neq 0$.}
Similarly as for the case $\mu =0$, we take different values for the parameters $k,\rho_{1},\rho_{2},\rho_{3},b,\beta $ and $\tau $. For example in Figure 3, we take $k=b=\rho_{1}=\rho_{2}=2$ and $\rho_{3}=\delta =\beta=\tau  =1$.
\begin{figure}[H]
        \centering
        \includegraphics[width=4cm]{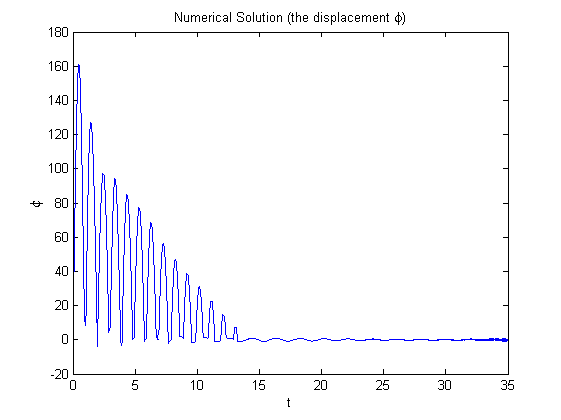}
        \caption{\label{fig3}}
\end{figure}
Here we have theoretically obtained a polynomial decay of the energy.\\
Finally, in figure 4 we give the three dimensional pointwise numerical solution of the Timoshenko system \eqref{1}. This proves again the energy decay of the transversal displacement $\varphi $, for $t$ large enough.
\begin{figure}[H]
        \includegraphics[width=8cm]{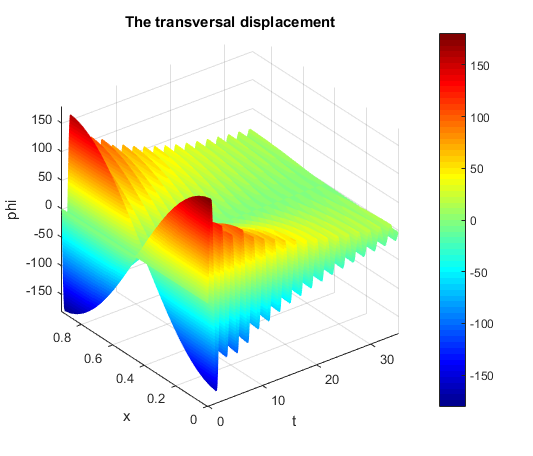}
        \caption{}
\end{figure}
%\textbf{Conclusion}
%
%$\bullet $ The function curve decreases more rapidly in Figure \ref{fig1} as in Figure \ref{fig3} which is proved theoretically.
%
%$\bullet $ We can choose other numerical test, but we must choose the constant c for convergence of the solution.
%%%--------------------

%%---------------------------

\begin{thebibliography}{99}

\bibitem{Ayadi} Ayadi M.A, Bchatnia A, Hamouda M and Messaoudi S.,
\textit{General decay in a Timoshenko-type system with thermoelasticity with second sound}. Advances Nonlinear Analysis, DOI: 10.1515/anona-2015-0038.
\bibitem{Raposo} C. A. Raposo, J. A. D. Chuquipoma, J. A. J. Avila, M. L. Santos, \textit{Exponential decay and numerical solution for a Timoshenko system with delay term in the internal feedback. International Journal of Analysis and Applications.} Vol. 3, no. \textbf{1}, (2013), 1-13.
\bibitem{FernandezRacke}   Fern\'{a}ndez Sare, H. D., and Racke, R.,  \textit{On the stability of damped
Timoshenko systems: Cattaneo versus Fourier law}, Arch. Rational Mech. Anal.,
\textbf{194} (\textbf{1}) (2009), 221-251.
\bibitem{SoufyaneWehbeUnifo}  Soufyane A. and Wehbe A., \textit{Uniform stabilization for the
Timoshenko beam by a locally distributed damping,} Electron. J. Differential
Equations no. \textbf{29} (2003), 1-14.
\bibitem{MessaoudiMustafa3}  Messaoudi, S. A., and Mustafa, M. I.,  \textit{On the stabilization of the
Timoshenko system by a weak nonlinear dissipation}, Math. Meth. Appl. Sci.,
\textbf{32} (\textbf{4}) (2009), 454-469.
\bibitem{MunozRiveraRackeMildy}  Mu\~{n}oz Rivera J.E. and Racke R., \textit{Mildly dissipative
nonlinear Timoshenko systems-global existence and exponential stability, J}.
Math. Anal. Appl. \textbf{276} (2002), 248-276.

\end{thebibliography}
\end{document}